\documentclass[12pt,draft]{amsart}
\usepackage{amssymb,times,epsfig}

\makeatletter
\def\@strippedMR{}
\def\@scanforMR#1#2#3\endscan{
  \ifx#1M\ifx#2R\def\@strippedMR{#3}
  \else\def\@strippedMR{#1#2#3}
  \fi\fi}
\renewcommand\MR[1]{\relax\ifhmode\unskip\spacefactor3000 \space\fi
  \@scanforMR#1\endscan
  MR\MRhref{\@strippedMR}{\@strippedMR}}
\makeatother

\newcounter{roman}

\newtheorem*{Thm*}{Theorem}
\newtheorem{Thm}{Theorem}
\newtheorem{Cor}[Thm]{Corollary}
\newtheorem{Prop}[Thm]{Proposition}
\newtheorem{Lemma}[Thm]{Lemma}

\theoremstyle{definition}
\newtheorem{Defn}{Definition}
\newtheorem{Notation}[Defn]{Notation}

\newtheorem{Remark}{Remark}
\newtheorem{Example}{Example}

\newcommand{\mf}[1]{\mathbb{#1}}
\newcommand{\mc}[1]{\mathcal{#1}}

\DeclareMathOperator{\tr}{\mathrm{tr}}

\newcommand{\norm}[1]{\left\Vert#1\right\Vert}
\newcommand{\abs}[1]{\left\vert#1\right\vert}
\newcommand{\chf}[1]{\mathbf{1}_{#1}}
\newcommand{\set}[1]{\left\{#1\right\}}

\newcommand{\State}[1]{\Phi \left[ #1 \right]}
\renewcommand{\phi}{\varphi}
\newcommand{\eps}{\varepsilon}

\renewcommand{\Re}{\mathrm{Re}\ }
\renewcommand{\Im}{\mathrm{Im}\ }
\newcommand{\supp}[1]{\mathrm{supp} \left( #1 \right)}

\newcommand{\br}{\medskip\noindent}

\allowdisplaybreaks[1]

\title{Bochner-Pearson-type characterization of the free Meixner class}
\author[M.~Anshelevich]{Michael Anshelevich}
\address{Department of Mathematics, Texas A\&M University, College Station, TX 77843-3368}
\email{manshel@math.tamu.edu}
\subjclass[2000]{Primary 46L54; Secondary 05E35, 33C45, 34B24}
\date{\today}

\begin{document}

\begin{abstract}
The operator $L_\mu: f \mapsto \int \frac{f(x) - f(y)}{x - y} \,d\mu(y)$ is, for a compactly supported measure $\mu$ with an $L^3$ density, a closed, densely defined operator on $L^2(\mu)$. We show that the operator $Q = p L_\mu^2 - q L_\mu$ has polynomial eigenfunctions if and only if $\mu$ is a free Meixner distribution. The only time $Q$ has orthogonal polynomial eigenfunctions is if $\mu$ is a semicircular distribution. More generally, the only time the operator $p (L_\nu L_\mu) - q L_\mu$ has orthogonal polynomial eigenfunctions is when $\mu$ and $\nu$ are related by a Jacobi shift.
\end{abstract}

\maketitle

\section{Introduction}

\br
The Meixner class of orthogonal polynomials, and its corresponding family of distributions---Gaussian, Poisson, gamma, binomial, negative binomial, and hyperbolic secant---has a number of combinatorial and probabilistic characterizations. For example, in the original paper \cite{Meixner}, Meixner showed that this class consists exactly of all orthogonal polynomials whose exponential generating functions have the special form $A(z) e^{x B(z)}$. For probabilistic characterizations, see for example \cite{Laha-Lukacs,Wes-commutative}. This class also plays an important role in statistics, under the name of quadratic exponential families, see \cite{Morris,Morris-Statistical,Diaconis-Gibbs-sampling}.

\br
On the other hand, historically, the most important class of orthogonal polynomials, often called the ``classical orthogonal polynomials'' \cite{Al-Salam-Characterization}, consists of the Hermite, Laguerre, and Jacobi families. This class also has numerous characterizations, notably the following one in terms of second order differential equations \cite[Chapter 20]{Ismail-Encyclopedia}.

\begin{Thm*}[Bochner 1929]
The Sturm-Liouville operator
\[
p(x) y'' + q(x) y'
\]
has a family of eigenfunctions consisting of orthogonal polynomials if and only if these eigenfunctions are Hermite, Laguerre, or Jacobi. In this case $p$ is a polynomial of degree at most $2$, and $q$ a polynomial of degree at most $1$.
\end{Thm*}

\br
In all of these cases, the operator is symmetric with respect to the inner product induced by the orthogonality measure of the polynomials. This orthogonality measure has a density determined by
\[
\frac{w'}{w} = \frac{q - p'}{p}.
\]
The class of probability densities with $\frac{w'}{w} = \frac{\text{linear}}{\text{quadratic}}$ is larger than the Bochner class, and is known in statistics as the Pearson class \cite{Diaconis-Closed-form-summation}; the corresponding operator has polynomial eigenfunctions, which however are not necessarily orthogonal on the real line. See Remark~\ref{Remark:Bochner} for more details. Note that all measures in this class have continuous densities, while the Meixner class includes purely atomic measures such as the Poisson distribution.

\br
With motivation from free probability \cite{VDN,Nica-Speicher-book}, in \cite{AnsMeixner} I introduced the free Meixner distributions as all the probability measures whose orthogonal polynomials have the ordinary generating function of the special form $A(z) \frac{1}{1 - x B(z)}$. The free Meixner class also has an explicit description, and is in fact quite simple, simpler than the classical Meixner class \cite{CTConstant}. Nevertheless, the majority of characterizations of the classical Meixner class have their counterparts---some straightforward, some not---for the free Meixner class \cite{Boz-Bryc,Bryc-Cauchy-Stieltjes,AnsFree-Meixner,Bozejko-Lytvynov-Meixner-I}. The goal of this paper is to investigate the free version of the Bochner / Pearson class. We will show that, in the appropriate framework, the Bochner-type characterization holds for the free Meixner class itself, quite unlike in the classical situation.

\br
What is the ``free version'' of the derivative operator? Unfortunately (or, rather, fortunately), there is no classical-free theory dictionary. However, Voiculescu (in \cite{Voi-Entropy5,VoiCoalgebra} and subsequent work) has shown that some roles played by the derivative operator are played in free probability by the difference quotient operator,
\[
\partial: C[\mf{R}] \rightarrow C[\mf{R} \times \mf{R}], \qquad (\partial f)(x,y) = \frac{f(x) - f(y)}{x-y}.
\]
It is hard to talk about differential equations in the difference quotient operator, since its domain and range are not the same. One modification of it is the single variable operator
\[
D_0 f(x) = \frac{f(x) - f(0)}{x},
\]
which is the $q=0$ version of the $q$-derivative operator from special function theory. The answer for the Bochner-type question for the $q$-derivative $D_q$ is known, see \cite{Al-Salam-Characterization} or \cite[Chapter 20]{Ismail-Encyclopedia}. However, for $q=0$, all these classes of polynomials reduce simply to monomials (with the exception of the degree one polynomial), and so are quite uninteresting. Other questions of this type involve the Askey-Wilson operator $\mc{D}_q$, and the difference operator. In the latter case, the corresponding orthogonal polynomials are Charlier, Meixner and Krawtchouk (that is, the discrete measures in the Meixner class), but also the Hahn polynomials. Note that the Meixner-Pollaczek polynomials belong to the Meixner class without satisfying either a differential or a (real) difference equation.

\br
Instead of the choices above, we will consider the following single variable operator: for a measure $\mu$,
\[
L_\mu[f] = (I \otimes \mu) [\partial f] = \int_{\mf{R}} \frac{f(x) - f(y)}{x-y} \,d\mu(y).
\]
So this is really an integral operator, but we will think of it as an analog of the derivative operator. Of course, the big difference from the derivative or the $D_0$ operator is that it depends on $\mu$. (Note, however, that $D_0 = L_{\delta_0}$.) Our Bochner-type question is: for which $\mu$ can a second order operator in $L_\mu$ with polynomial coefficients have polynomial eigenfunctions? And when are these eigenfunctions orthogonal?

\br
Here is the summary of the main results of this paper. Section~\ref{Section:Preliminaries} contains, in particular, the detailed description of the Bochner and Pearson classes on one hand, and of the free Meixner class on the other, as a convenient reference. In Section~\ref{Section:free-Meixner}, we show that an operator
\[
p L_\mu^2 + q L_\mu
\]
has polynomial eigenfunctions if and only if $\mu$ has a free Meixner distribution (given $\mu$, the choice of $p, q$ is unique). On the other hand, if $p = 1$ and $q = - H$, the conjugate variable for $\mu$, then the operator $p L_\mu^2 + q L_\mu$ is symmetric. $H$ is a polynomial of degree one exactly when $\mu$ is a semicircular distribution, which are also the only cases when the polynomial eigenfunctions are orthogonal with respect to $\mu$. Finally, in Section~\ref{Section:Extensions} we look at operators of the form $p (L_\nu L_\mu) + q L_\mu$ which depend on a \emph{pair} of measures. We classify pairs $(\mu, \nu)$ which correspond to polynomial eigenfunctions. For any $\mu$, if $\nu$ is a Bernoulli distribution,  one can choose the coefficients in the operator so that it will have polynomial eigenfunctions. On the other hand, for any $\nu$, if $\mu$ is a ``Jacobi shift'' of $\nu$, the operator has polynomial eigenfunctions. Finally, we show that the eigenfunctions are orthogonal with respect to $\mu$ if and only if the pair $(\mu, \nu)$ is of the latter type; in that case, the operator is also symmetric.

\br
\textbf{Acknowledgements.} Thanks to Alberto Gr{\"u}nbaum for questions that motivated the writing of this paper, and to Serban Belinschi and David Damanik for useful conversations. Thanks also to Dennis Stanton, the person who knows everything about orthogonal polynomials, new and old, and generously shares it.

\section{Preliminaries}
\label{Section:Preliminaries}

\begin{Notation}
A \emph{polynomial system} is a family of polynomials $\set{P_n: n \geq 0}$ with $\deg P_n = n$.
\end{Notation}

\noindent
We start with a more detailed description of Bochner's theorem, to be compared with the results that follow it.

\begin{Remark}[Bochner-Pearson classification]
\label{Remark:Bochner}
If the operator
\[
y \mapsto (p(x) y')' + q(x) y'
\]
has a polynomial system of eigenfunctions, then $p$ is a polynomial of degree at most $2$, and $q$ a polynomial of degree at most $1$, and the operator is symmetric with respect to the inner product induced by the measure $d\mu(x) = w(x) \,dx$ with the density determined by
\[
\frac{w'}{w} = \frac{q}{p}.
\]
By affine transformations, one can reduce the analysis to the following cases.
\begin{enumerate}
\item
$\frac{w'}{w} = - \alpha x$. Then $\alpha > 0$ (and we may taken $\alpha = 1$), and the Hermite polynomials are eigenfunctions, orthogonal with respect  to the Gaussian distribution $\frac{1}{2 \pi} e^{-x^2/2}$.
\item
$\frac{w'}{w} = \frac{\alpha - 1}{x} - \beta$. Then $\beta > 0$ (and can be taken $\beta = 1$), $\alpha > 0$, and the Laguerre polynomials are eigenfunctions, orthogonal with respect to the gamma distribution $\frac{1}{\Gamma(\alpha)} x^{\alpha - 1} e^{-x}$ on $(0, \infty)$.
\item
$\frac{w'}{w} = \frac{\alpha - 1}{x} - \frac{\beta - 1}{1 - x}$. Then $\alpha, \beta > 0$, and the Jacobi polynomials are eigenfunctions, orthogonal with respect to the beta distribution $\frac{\Gamma(\alpha + \beta)}{\Gamma(\alpha) \Gamma(\beta)} x^{\alpha - 1} (1 - x)^{\beta - 1}$ on $(0,1)$.
\item
$\frac{w'}{w} = - \frac{1 + \alpha}{x} + \frac{\beta}{x^2}$, $\alpha > 0$, $\beta > 0$ so take $\beta = 1$, and the measure is $\frac{1}{Z} x^{- 1 - \alpha} e^{-1/x}$ on $(0, \infty)$ ($Z$ will always denote the appropriate normalization constant). Bessel polynomials, orthogonal in the complex plane, are eigenfunctions.
\item
$\frac{w'}{w} = - \frac{(1 + \alpha) x}{x^2 + 1} + \frac{\beta}{x^2 + 1}$, $\alpha > 0$, and the measure is $\frac{1}{Z} (x^2 + 1)^{- (1 + \alpha)/2} e^{\beta \tan^{-1}(x)}$. For $\beta = 0$, these are the $t$-distributions.
\end{enumerate}
\end{Remark}

\begin{Remark}
While there is, in general, no relation between the Pearson and Meixner distributions, some parallels between these classes can be explained by the following observation. As mentioned above, for the Pearson class $d\mu(x) = w(x) \,dx$,
\[
\frac{w'(x)}{w(x)} = \frac{d + e x}{1 + b x + c x^2}.
\]
On the other hand, let $\set{P_n(x,t)}$ be the monic orthogonal polynomials for the measure $\mu$ ($t$ is the convolution parameter), and denote by $F$ their exponential generating function:
\[
F(x,t,z) = \sum_{n=0}^\infty \frac{1}{n!} P_n(x,t) z^n.
\]
$\mu$ is in the Meixner class if and only if
\[
\frac{\partial_z F}{F} = \frac{x - t z}{1 + b z + c z^2}.
\]
\end{Remark}

\begin{Notation}
For a probability measure $\mu$, its Cauchy transform is
\[
G_\mu(z) = \int_{\mf{R}} \frac{1}{z - x} \,d\mu(x).
\]
Recall that $\mu$ can be recovered from it as the weak limit
\[
d\mu(x) = - \frac{1}{\pi} \lim_{\eps \rightarrow 0^+} \Im G(x + i \eps) \,dx.
\]
The Hilbert transform of $\mu$, defined Lebesgue almost everywhere, is the function
\[
H[\mu](x) = \frac{1}{\pi} \lim_{\eps \rightarrow 0^+} \Re G(x + i \eps).
\]
If $\mu$ has a density in $L^p(\mf{R}, dx)$, $p > 1$, then $H[\mu] \in L^p(\mf{R}, dx)$ \cite[Chapter~4]{Tricomi-Integral-equations}.

\br
The $R$-transform of $\mu$ is determined by
\[
G \left( R(z) + \frac{1}{z} \right) = z;
\]
it is an analytic function on a domain, but in this paper we will identify it with its power series expansion.

\br
Throughout most of the paper, $\mu$ (and later, also $\nu$) will be a compactly supported probability measure. In particular, its moments
\[
m_n = \mu[x^n]
\]
are finite, and we will also consider $\mu$ as a linear functional on polynomials. In this case, its moment generating function is
\[
M(z) = M_\mu(z) = \sum_{n=0}^\infty m_n z^n.
\]
\end{Notation}

\begin{Remark}
\label{Remark:Conjugate}
For a compactly supported probability measure $\mu$ with an $L^3(\mf{R}, dx)$ density, its renormalized Hilbert transform
\[
H = H_\mu = 2 \pi H[\mu]
\]
is well defined and is in $L^2(\mu)$. Moreover, for $f \in L^2(\mu)$,
\[
(\mu \otimes \mu) [\partial f] = \mu[H f].
\]
In other words, $H$  is the conjugate variable for $\mu$; see \cite{Voi-Entropy5} for more details. As explained there, the conjugate variable is the free analog of the classical score function $-\frac{w'}{w}$ of a probability measure $w(x) \,dx$.
\end{Remark}

\begin{Remark}[Free Meixner distributions]
\label{Remark:Meixner}
The normalized free Meixner distributions $\mu_{b,c}$ are probability measures with Jacobi parameter sequences
\[
\set{(0, b, b, b, \ldots), (1, 1+c, 1+c, 1+c, \ldots)},
\]
for $b \in \mf{R}$, $c \geq -1$. The general free Meixner distributions are affine transformations of these. More explicitly, the distribution with parameters $b, c$ is
\[
\frac{1}{2 \pi} \cdot \frac{\sqrt{\Bigl(4 (1 + c) - (x - b)^2\Bigr)_+}}{1 + b x + c x^2} \,dx + 0, 1, \text{or } 2 \text{ atoms}.
\]
Unfortunately, none of the descriptions of these measures in \cite{CTConstant,SaiConstant,Bryc-Cauchy-Stieltjes}, including our own \cite{AnsMeixner}, are complete, so we provide a detailed description here. By affine transformations, the situation can be reduced to the following six cases. For future reference, we also record their conjugate variables.
\begin{list}
{\textbf{\roman{roman}.}}{\usecounter{roman}}
\item
\label{semicircular}
$b = c = 0$,
\[
d\mu(x) = \frac{1}{2 \pi} \sqrt{(4 - x^2)_+} \,dx.
\]
This is the semicircular distribution. Conjugate variable $x$.
\item
\label{Marchenko}
$c = 0$, $b \neq 0$, $\alpha > 0$,
\[
\frac{1}{2 \pi} \frac{\sqrt{\Bigl((1 + \sqrt{\alpha})^2 - x \Bigr) \Bigl(x - (1 - \sqrt{\alpha})^2 \Bigr)_+}}{x} \,dx + \max(1 - \alpha, 0) \delta_0.
\]
These are the Marchenko-Pastur distributions. Conjugate variable $\frac{1 - \alpha}{x} + 1$.
\item
\label{binomial}
$-1 \leq c < 0$, $\alpha, \beta > 0$, $\alpha + \beta \geq 1$,
\[
\frac{1}{2 \pi} \frac{\sqrt{\Bigl(4 \alpha (1 - x) - (\alpha + 1 - (\alpha + \beta) x)^2 \Bigr)_+}}{x (1 - x)} \,dx + \max(1 - \alpha, 0) \delta_0 + \max(1 - \beta, 0) \delta_1.
\]
These are the free binomial distributions, including the Bernoulli distributions for $c = -1$, $\alpha + \beta = 1$. Conjugate variable $\frac{1 - \alpha}{x} - \frac{1 - \beta}{1 - x}$.
\item
\label{gamma}
$c \neq 0$, $b^2 - 4 c = 0$, $\alpha > 0$,
\[
\frac{1}{2 \pi} \frac{\sqrt{\Bigl( (1 + \sqrt{\alpha + 1})^2 - \alpha x \Bigr) \Bigl( \alpha x - (1 - \sqrt{\alpha + 1})^2 \Bigr)_+}}{x^2} \,dx.
\]
Conjugate variable $\frac{2 + \alpha}{x} - \frac{1}{x^2}$.
\item
\label{secant}
$c > 0$, $b^2 - 4 c < 0$, $\alpha > 0$, $\beta \in \mf{R}$,
\[
\frac{1}{2 \pi} \frac{\sqrt{\Bigl( \alpha^2 x^2 - 2 \beta (2 + \alpha) x + 4 (\beta^2 - 1 - \alpha) \Bigr)_+}}{1 + x^2} \,dx.
\]
Conjugate variable $\frac{(2 + \alpha) x}{1 + x^2} - \frac{\beta}{1 + x^2}$.
\item
\label{negative}
$c > 0$, $b^2 - 4 c > 0$, $\alpha + \beta < 0$, $\alpha \beta < 0$,
\[
\begin{split}
\frac{1}{2 \pi} & \frac{\sqrt{\Bigl(4 \alpha (1 - x) - (\alpha + 1 - (\alpha + \beta) x)^2 \Bigr)_+}}{x (x - 1)} \,dx \\
& + \Bigl( (1 - \alpha) \delta_0 \text{ for } 0 < \alpha < 1 \Bigr) + \Bigl( (1 - \beta) \delta_1 \text{ for } 0 < \beta < 1 \Bigr);
\end{split}
\]
note that since $\alpha + \beta < 0$, at most one atom may appear. Conjugate variable $\frac{1 - \alpha}{x} - \frac{1 - \beta}{1 - x}$.
\end{list}

\br
Using the $(b,c)$ parametrization, and either the orthogonal polynomials \cite{AnsFree-Meixner} or the Laha-Lukacs \cite{Boz-Bryc} characterization, the corresponding classical Meixner classes are (\textbf{i}) Gaussian (\textbf{ii}) Poisson (\textbf{iii}) binomial (\textbf{iv}) gamma (\textbf{v}) hyperbolic secant (\textbf{vi}) negative binomial. The only difference is that in the binomial case, $- \frac{1}{c} \in \mf{N}$, but for the free binomial one can take any $-1 \leq c < 0$. See also Remark~\ref{Remark:Random-matrices}.

\br
Another description of the free Meixner distributions \cite{Boz-Bryc,AnsFree-Meixner}, again with a slightly different normalization, is in terms of the $R$-transform: for the distribution with parameters $b, c$, mean $m$ and variance $t$,
\begin{equation}
\label{Characteristic-equation}
\frac{R(u) - m}{u} = t + b (R(u) - m) + (c/t) (R(u) - m)^2.
\end{equation}
Note that the second paper quoted above uses the combinatorial $R$-transform $\mc{R}(u) = u R(u)$.
\end{Remark}

\begin{Remark}
Our main object it the operator
\[
L_\mu[f] = (I \otimes \mu)[\partial f] = \int_{\mf{R}} \frac{f(x) - f(y)}{x - y} \,d\mu(y).
\]
We will consider three versions of it. First, $L_\mu$ as an operator on the vector space of polynomials,
\[
L_\mu: \mf{R}[x] \rightarrow \mf{R}[x], \qquad L_\mu[x^n] = \sum_{k=0}^{n-1} m_{n-k-1} x^k.
\]
Like other deformations of the derivative operator, it lowers the degree by $1$. See also Lemma~\ref{Lemma:free-Appell}.

\br
Second, $L_\mu$ is a linear operator on the vector space of continuously differentiable functions, mapping it to continuous functions,
\[
L_\mu: C^{(1)}(\mf{R}) \rightarrow C(\mf{R}).
\]
On this space, it has eigenfunctions
\[
L_\mu \left[ \frac{1}{z - x} \right] = \mu\left[ \frac{1}{z - x} \right] \frac{1}{z - x} = G_\mu(z) \frac{1}{z - x}.
\]
For real $z$, the eigenvalue is $\pi H[\mu](z)$, where $H$ is the Hilbert transform, whenever this quantity is finite.

\br
The third version is described in the following proposition. See also Lemma~\ref{Lemma:Self-adjoint}.
\end{Remark}

\begin{Prop}
\label{Prop:L-closed}
Let $\mu$ be a compactly supported measure with density $w \in L^3(\mf{R}, dx)$. The operator $L_\mu$ is a bounded operator from $L^2(\mu)$ to $L^1(\mu)$ and from $L^\infty(\mu)$ to $L^2(\mu)$. It is a (possibly) unbounded, densely defined, closed linear operator on $L^2(\mu)$.
\end{Prop}

\begin{proof}
Denoting $\norm{H}_{L^p(\mf{R}, dx) \rightarrow L^p(\mf{R}, dx)} = C_p$,
\[
\mu[H[\mu]^2] = \norm{H[\mu]^2 w}_1 \leq \norm{H[\mu]}_3^2 \norm{w}_3 \leq C_3^2 \norm{w}_3^3.
\]
So for $f \in L^2(\mu)$,
\[
\mu[\abs{f H[\mu]}] \leq \sqrt{\mu[f^2] \mu[H[\mu]^2]} \leq C_3 \norm{w}_3^{3/2} \mu[f^2]^{1/2}.
\]
Also,
\[
\mu[\abs{H[f \mu]}] = \norm{H[f \mu] w}_1 \leq \norm{H[f \mu]}_{3/2} \norm{w}_3 \leq C_{3/2} \norm{f w}_{3/2} \norm{w}_3.
\]
But
\[
\begin{split}
\norm{f w}_{3/2} = \norm{\abs{f}^{3/2} \abs{w}^{3/2}}_1^{2/3} & \leq \norm{\abs{f}^{3/2} w^{3/4}}_{4/3}^{2/3} \norm{w^{3/4}}_4^{2/3} \\
& = \norm{f^2 w}^{1/2} \norm{w^3}^{1/6} = \mu[f^2]^{1/2} \norm{w}_3^{1/2}.
\end{split}
\]
Thus
\[
\mu[\abs{H[f \mu]}] \leq C_{3/2} \norm{w}_3^{3/2} \mu[f^2]^{1/2}.
\]
It follows that all the quantities in the following equation are well defined,
\[
L_\mu[f] = \pi f(x) H[\mu](x) - \pi H[f \mu](x),
\]
and $\norm{L_\mu}_{L^2(\mu) \rightarrow L^1(\mu)} \leq \pi (C_3 + C_{3/2}) \norm{w}_3^{3/2}$. Similarly, for $f \in L^\infty(\mu)$,
\[
\mu[H[f \mu]^2]^{1/2} = \norm{H[f \mu]^2 w}_1^{1/2} \leq \norm{H[f \mu]^2}_{3/2}^{1/2} \norm{w}_3^{1/2}
= \norm{H[f \mu]}_3 \norm{w}_3^{1/2} \leq C_3 \norm{f}_\infty \norm{w}_3^{3/2}
\]
and $\norm{L_\mu}_{L^\infty(\mu) \rightarrow L^2(\mu)} \leq 2 \pi C_3 \norm{w}_3^{3/2}$.

\br
Finally, let $f_n \rightarrow f$ in $L^2(\mu)$ and $g_n = L_\mu[f_n] \rightarrow g$ in $L^2(\mu)$. Then $g_n \rightarrow L_\mu[f]$ in $L^1(\mu)$, so $g = L_\mu[f]$, and $L_\mu : L^2(\mu) \rightarrow L^2(\mu)$ is closed.
\end{proof}

\section{The free Meixner characterization}
\label{Section:free-Meixner}

\begin{Thm}
\label{Thm:SL-one-measure}
For $\mu$ with all moments finite, consider a ``Sturm-Liouville''-type operator of the form
\[
Q_\mu = p(x) L_\mu^2 + q(x) L_\mu,
\]
where $p(x), q(x)$ are polynomials. It has a polynomial system of eigenfunctions only if $\mu$ is a free Meixner distribution. In that case, the conjugate variable for $\mu$ is $- \frac{q(x)}{p(x)}$.
\end{Thm}

\begin{proof}
Clearly
\[
Q[x] = q(x).
\]
Thus to have a polynomial eigenfunction of degree $1$, we need $q$ to have degree at most $1$. This implies that $p$ has degree at most $2$. Thus denote
\[
p(x) = a + b x + c x^2
\]
and
\[
q(x) = d + e x.
\]
Next, we compute
\[
Q[1] = 0,
\]
\[
Q[x] = (d + e x)
\]
and
\[
Q[x^n] = (a + b x + c x^2) (x^{n-2} + \ldots) + (d + e x) (x^{n-1} + \ldots) = (c + e) x^n + \ldots.
\]
Thus the eigenvalues (if any) are $0$ for $n=0$, $e$ for $n=1$, and $c + e$ for $n \geq 2$. Therefore the polynomial system of eigenfunctions can be taken to be
\[
\set{1, x + \alpha, x^n + \beta_n x + \gamma_n, n \geq 2}
\]
for some $\alpha, \beta_n, \gamma_n$. We compute
\[
Q[x + \alpha] = d + e x = e (x + \alpha),
\]
thus
\begin{equation}
\label{Linear-coefficient}
e \alpha = d.
\end{equation}
Similarly, for $n \geq 2$
\[
\begin{split}
Q[x^n + \beta_n x + \gamma_n]
& = (a + b x + c x^2) \sum_{k=0}^{n-2} x^k \sum_{i=0}^{n-k-2} m_i m_{n-k-i-2} + (d + e x) \left( \sum_{k=0}^{n-1} x^k m_{n-k-1} + \beta_n \right) \\
& = (c + e) (x^n + \beta_n x + \gamma_n).
\end{split}
\]
Comparing coefficients, we get for $k=0$
\begin{equation}
\label{Gamma}
a \sum_{i=0}^{n-2} m_i m_{n-i-2} + d m_{n-1} + d \beta_n = (c + e) \gamma_n,
\end{equation}
for $k=1$
\begin{equation}
\label{Beta}
a \sum_{i=0}^{n-3} m_i m_{n-i-3} + b \sum_{i=0}^{n-2} m_i m_{n-i-2} + d m_{n-2} + e m_{n-1} + e \beta_n = (c + e) \beta_n,
\end{equation}
while for $2 \leq k \leq n-1$
\[
a \sum_{i=0}^{n-k-2} m_i m_{n-k-i-2} + b \sum_{i=0}^{n-k-1} m_i m_{n-k-i-1} + c \sum_{i=0}^{n-k} m_i m_{n-k-i} + d m_{n-k-1} + e m_{n-k} = 0.
\]
Here the empty sums are understood to be zero. Changing variables from $n-k$ to $n$, we see that for $n \geq 1$
\begin{equation}
\label{Moments-recursion}
a \sum_{i=0}^{n-2} m_i m_{n-i-2} + b \sum_{i=0}^{n-1} m_i m_{n-i-1} + c \sum_{i=0}^{n} m_i m_{n-i} + d m_{n-1} + e m_{n} = 0.
\end{equation}
Note that for $n=0$, the corresponding term is $c + e$. In terms of the moment generating function $M(z)$ of $\mu$, equation~\eqref{Moments-recursion} gives
\[
a z^2 M(z)^2 + b z M(z)^2 + c M(z)^2 + d z M(z) + e M(z) = c + e.
\]
In terms of the Cauchy transform
\[
G_\mu(z) = \int_{\mf{R}} \frac{1}{z - x} \,d\mu(x) = \frac{1}{z} M(1/z),
\]
(treated as a formal power series if necessary) the relation is
\begin{equation}
\label{Cauchy}
(a + b z + c z^2) G(z)^2 + (d + e z) G(z) - (c + e) = 0.
\end{equation}
So
\[
G(z) = \frac{- (d + e z) - \sqrt{(d + e z)^2 + 4 (c + e) (a + b z + c z^2)}}{2 (a + b z + c z^2)}.
\]
This means that on the support of $\mu$
\[
\set{x: (d + e x)^2 + 4 (c + e) (a + b x + c x^2) \leq 0} \cup \set{x: a + b x + c x^2 = 0},
\]
the conjugate variable of $\mu$ is
\[
- \frac{d + e x}{a + b x + c x^2}.
\]
On the other hand, in terms of the $R$-transform determined by
\[
G \left( R(z) + \frac{1}{z} \right) = z,
\]
equation~\eqref{Cauchy} states that
\[
c u R(u)^2 + (b u + (2c + e)) R(u) + (a u + b + d) = 0,
\]
or
\begin{equation}
\label{Riccati}
-(2 c + e) \frac{R(u)}{u} - \frac{b + d}{u} = a + b R(u) + c R(u)^2,
\end{equation}
which is equation~\eqref{Characteristic-equation}.

\br
For future reference, we also compute the eigenfunctions. Denoting
\[
B(z) = \sum_{n=2}^\infty \beta_n z^n,
\]
we get from equation~\eqref{Beta} that
\[
c B(z) = a z^3 M(z)^2 + b z^2 M(z)^2 + d z^2 M(z) + e z (M(z) - 1) = (c + e) z - c z M(z)^2 - e z,
\]
so either $c=0$ or
\[
B(z) = z (1 - M(z)^2).
\]
Similarly, for
\[
C(z) = \sum_{n=2}^\infty \gamma_n z^n,
\]
we get
\[
(c + e) C(z) = a z^2 M(z)^2 + d z (M(z) - 1) + d B(z) =  z M(z) (a z M(z) + d(1 - M(z)))
\]
if $c \neq 0$.
\end{proof}

\begin{Prop}
If $\mu$ is a normalized free Meixner distribution with parameters $b$ and $c \geq -1$, the corresponding operator $Q$ from Theorem~\ref{Thm:SL-one-measure} has polynomial eigenfunctions for the unique (up to a factor) choice of
\[
p(x) = 1 + b x + c x^2, \qquad q(x) = - (b + (1 + 2 c) x).
\]
If $c = - \frac{1}{2}$, we additionally require that $b = 0$.
\end{Prop}

\begin{proof}
By shifting and re-scaling the measure, we may assume that $m_1 = 0$, $m_2 = 1$. From equation~\eqref{Moments-recursion},
\[
b + 2 m_1 c + d + m_1 e = 0
\]
and
\[
a + 2 m_1 b + (2 m_2 + m_1^2) c + m_1 d + m_2 e = 0.
\]
So under this normalization,
\[
b + d = 0, \qquad a + 2 c + e = 0.
\]
In equation \eqref{Riccati}, this corresponds to
\[
a \frac{R(u)}{u} = a + b R(u) + c R(u)^2.
\]
If $a = 0$, then $R(u) (b + c R(u)) = 0$, which corresponds to constant $R(u)$ and $\mu$ being a delta measure. If $a \neq 0$, we may assume without loss of generality that $a > 0$. In this case
\[
\frac{R(u)}{u} = 1 + (b/a) R(u) + (c/a) R(u)^2,
\]
which is equation~\eqref{Characteristic-equation} in standard form (for mean zero and variance one). In particular, we know that $R(u)$ is an $R$-transform of a positive measure for any $b$ and for $c/a \geq -1$, in other words for
\[
a + c \geq 0.
\]
Note that this implies the eigenvalue $c + e = - (a + c) \leq 0$, while $e = - (a + 2 c)$ may be positive, negative, or zero.

\br
For $e = - (a + 2 c) = 0$, equation~\eqref{Linear-coefficient} additionally implies that $b = -d = 0$ as well.
\end{proof}

\begin{Example}
We thus get a correspondence between operators
\[
(p y')' + q y'
\]
and
\[
p L_\mu^2 + q L_\mu
\]
having polynomial eigenvalues, given by
\[
\frac{q}{p} = \frac{w'}{w} = - H_\mu.
\]
Note that the operators have slightly different forms, but this is the appropriate correspondence. Note also that the parameter ranges in parameterizations in Remarks~\ref{Remark:Bochner} and \ref{Remark:Meixner} coincide in some but not all cases, and that classically, there is no analog of the restriction at the very end of the proof above.

\br
The special classes of measures are: semicircular for
\[
Q = L_\mu^2 - x L_\mu
\]
with eigenvalue $-1$, Marchenko-Pastur for
\[
Q = (a + b x) L_\mu^2 - (b + a x) L_\mu
\]
with eigenvalue $-a$, Bernoulli for
\[
Q = (a + b x - a x^2) L_\mu^2 - (b - a x) L_\mu
\]
with eigenvalues $a, 0$ and for general $a + c \geq 0$,
\[
Q = (a + b x + c x^2) L_\mu^2 - (b + (a + 2 c) x) L_\mu.
\]
If $a + 2 c = 0$, or in other words $(c/a) = -\frac{1}{2}$, then also $b=0$ from the condition $e \alpha = d$. The corresponding distribution is the arcsine law, with
\[
Q = (2 - x^2) L_\mu^2
\]
and eigenvalues $0, - \frac{1}{2}$. More generally, in category (\textbf{iii}) or Remark~\ref{Remark:Meixner}, the only values of $\alpha, \beta$ with  $\alpha + \beta = 2$ producing polynomial eigenfunctions are $\alpha = \beta = 1$. Note that classically, these parameter values correspond to the uniform distribution, and Legendre polynomials.
\end{Example}

\begin{Remark}
A very different characterization of the free Meixner class in terms of certain operators mapping polynomials to themselves appears in \cite{Bozejko-Lytvynov-Meixner-I}. It would be interesting to see if there is a relation to our results. Note however that the classical version of the results of that paper involves the Meixner and not the Bochner class.
\end{Remark}

\begin{Lemma}
The only operators of the type in Theorem~\ref{Thm:SL-one-measure} with polynomial eigenfunctions orthogonal with respect to $\mu$ correspond to the semicircular distributions.
\end{Lemma}

\begin{proof}
By re-scaling, we may assume that $\mu$ has mean zero and variance one. If the eigenfunction polynomials are orthogonal with respect to $\mu$, in particular they are centered with respect to it. So
\[
\mu[x + \alpha] = m_1 + \alpha = m_1 + (d/e) = 0
\]
implies $\alpha = 0$, and from \eqref{Linear-coefficient}, $0 = d = -b$. Also
\[
\mu[x^2 + \beta_2 x + \gamma_2] = m_2 + m_1 \beta_2 + \gamma_2 = m_2 - 2 (c/c) m_1^2 + \frac{a - d m_1}{c+e} = 0,
\]
so $a = - (c+e)$ and $c = 0$. Thus $\mu$ is the semicircular distribution.
\end{proof}

\begin{Remark}[Free probability]
\label{Remark:Free-probability}
Out of the five Pearson classes, only the first three correspond to orthogonal polynomials and the Bochner scheme. Interestingly, exactly the three classes corresponding to these under Theorem~\ref{Thm:SL-one-measure}, and not the rest of the free Meixner distributions, have an interpretation in free probability theory. The semicircle law is the ``free Gaussian'', since it appears as the limit in the free central limit theorem. Similarly, the Marchenko-Pastur law appears as the limit distribution in the free version of the Poisson limit theorem. Finally, the free binomial distribution is the sum of freely independent Bernoulli variables. See \cite{VDN} for more details.
\end{Remark}

\begin{Remark}[Random matrices]
\label{Remark:Random-matrices}
Let $V$ be a potential such that
\begin{equation}
\label{convex}
V \in C^{(2)}(\mf{R}), \quad V \text{ is convex,} \quad \int_{\mf{R}} e^{-V(x)} \,dx < \infty.
\end{equation}
Let $X$ be an Hermitian $n \times n$ random matrix \cite{Deift} distributed according to
\begin{equation}
\label{Potential}
\frac{1}{Z} e^{- \tr V(X)} \,dX.
\end{equation}
Then as $n \rightarrow \infty$, the spectral distribution of $X$ converges to a compactly supported measure $\mu$ such that
\begin{equation}
\label{Correspondence}
V' = 2 \pi H[\mu].
\end{equation}
So if $w(x) \,dx = e^{-V(x)} \,dx$, then $- \frac{w'}{w} = H$, the conjugate variable for $\mu$. Of course, this was one of the original motivations for Voiculescu's definition of the conjugate variable.

\br
As shown above, this is also exactly the correspondence between the Pearson and the free Meixner classes, although the relevant ranges of parameters do not coincide in all cases. $V$ with $V'(x) = \frac{d + e x}{a + b x + c x^2}$ satisfies condition \eqref{convex} only if it is an affine transformation of $V = x^2/2$, and if $X$ is in the Gaussian unitary ensemble
\[
\frac{1}{Z} e^{- \tr (X^2/2)} \,dX,
\]
then its spectral distribution converges to the semicircular distribution. But also, for $\alpha \geq 1$, $- (\alpha - 1) \log x + x$ satisfies \eqref{convex} on $(0, \infty)$, and if $X$ is in the Wishart ensemble \cite{Muirhead-book}
\[
\frac{1}{Z} \det X^{\alpha - 1} e^{- \tr (X)} \chf{0 \leq X} \,dX
\]
(note that in the $n=1$ case, this is the gamma distribution), then its spectral distribution converges to the Marchenko-Pastur law with parameter $\alpha$. Similarly, for $\alpha, \beta \geq 1$, $(\alpha - 1) \log x + (\beta - 1) \log (x-1)$ satisfies \eqref{convex} on $(0, 1)$, and if $X$ is in the Jacobi ensemble \cite{Capitaine-Casalis-Beta,Collins-Jacobi}
\[
\frac{1}{Z} \det (I - X)^{\alpha - 1} \det X^{\beta - 1} \chf{0 \leq X \leq I} \,dX
\]
(which is a matrix version of the beta distribution), the spectral distribution converges to a free binomial law with parameters $\alpha, \beta$. Note that these parameter values correspond precisely to the distributions being absolutely continuous. Finally, the remaining two distributions in the Pearson class do not correspond to convex potentials, but still satisfy the hypothesis of Theorem~1.3 of \cite[Chapter I]{Saff-Totik}, and so have a unique equilibrium measure, which however is not necessarily determined by equation~\eqref{Correspondence}.

\br
For $\alpha, \beta < 1$, one also has matrix models, but they are less canonical. Namely, let $P$ be an $n \times n$ projection matrix of rank $k \leq n$, and $X$ be a GUE matrix as above. Then the matrix
\[
W = X P X
\]
has the Wishart distribution with parameters $n, k$, and as $n \rightarrow \infty$ while $\frac{k}{n} \rightarrow \alpha$, the spectral distribution of $W$ converges to the Marchenko-Pastur distribution with parameter $\alpha$, $0 \leq \alpha \leq 1$ \cite{Muirhead-book}. Moreover, if $W_1, W_2$ are Wishart matrices with parameters $n, k_1$ and $n, k_2$, respectively, with $k_1, k_2 \leq n$, $k_1 + k_2 \geq n$, then
\[
(W_1 + W_2)^{-1/2} W_1 (W_1 + W_2)^{-1/2}
\]
is well defined, and as $n \rightarrow \infty$, $\frac{k_1}{n} \rightarrow \alpha$, $\frac{k_2}{n} \rightarrow \beta$, its spectral distribution converges to the free binomial distribution with parameters $0 \leq \alpha, \beta \leq 1$, $\alpha + \beta \geq 1$ \cite{Capitaine-Casalis-Beta}. There is also another model for this distribution involving a product of random projections \cite{Collins-Jacobi}.

\br
While there exist abstract characterizations of, for example, the Wishart distribution \cite{Bobecka-Wesolowski}, we are not aware of a Bochner type characterization for these matrix distributions for $2 \leq n < \infty$.
\end{Remark}

\br
We close the section with the following analog of the operator $i (y' + \frac{1}{2} \frac{w'}{w})$ being self adjoint with respect to $w(x) \,dx$; note the absence of boundary conditions.

\begin{Lemma}
\label{Lemma:Self-adjoint}
For $\mu$ with density in $L^3(\mf{R}, dx)$, the operator $A = i (L_\mu - \frac{1}{2} H[\mu])$ is a self-adjoint operator on $L^2(\mu)$ with dense domain $\mc{D} = \set{f \in L^2(\mu): A[f] \in L^2(\mu)}$.
\end{Lemma}

\begin{proof}
By Proposition~\ref{Prop:L-closed}, the operator $A$ maps $L^\infty(\mu)$ to $L^2(\mu)$, so the domain $\mc{D}$ is dense. Suppose that there are $g, h \in L^2(\mu)$ such that for all $f \in \mc{D}$,
\[
\mu[A[f] g] = \mu[f h].
\]
For all $f \in L^\infty(\mu)$,
\begin{equation*}
\begin{split}
\mu[L[f] g]
& = (\mu \otimes \mu) [(g \otimes I) \partial f]
= (\mu \otimes \mu)[\partial(gf) - (I \otimes f) \partial g] \\
& = \mu[H g f] - (\mu \otimes \mu)[(I \otimes f) \partial g] \\
& = \mu[H g f] - (\mu \otimes \mu)[(f \otimes I) \partial g]
= \mu[H g f] - \mu[f L[g]].
\end{split}
\end{equation*}
This implies that the operator $A$ is symmetric and, for such $f$,
\[
\mu[f h] = \mu[A[f] g] = \mu[f A[g]].
\]
Here we have also used the fact that $A$ maps $L^2(\mu)$ to $L^1(\mu)$. It follows that $h = A[g]$ in $L^1(\mu)$, so $A[g] \in L^2(\mu)$ and $g \in \mc{D}$. Therefore $A$ is self-adjoint.
\end{proof}

\section{Extensions}
\label{Section:Extensions}

\subsection{Higher-order operators.}
The operator
\[
Q = \sum_{k=0}^n p_k(x) L_\mu^k
\]
with
\[
p_k(x) = a_k x^k + \ldots
\]
has polynomial eigenfunctions if and only if
\[
\sum_{k=0}^n p_k(z) G_\mu(z)^k = \sum_{k=0}^n a_k.
\]
However, unlike in the second order case, where the complete description of free Meixner distributions is available, it is not clear which solutions of this equation are Cauchy transforms of positive measures. Such analysis is related to the study of free exponential families, see \cite{Bryc-Cauchy-Stieltjes}.

\subsection{The two measure case.}

\begin{Remark}
In \cite{AnsAppell3}, with motivation from the two-state free probability theory (also called the c-free theory), we introduced the c-free conjugate variable $H_{\mu, \nu}$: for a pair of measures $\mu, \nu$, $H_{\mu, \nu}$ (which may not exist) is determined by
\[
(\mu \otimes \nu)[\partial f] = \mu[H_{\mu, \nu} f].
\]
Explicitly, under appropriate conditions on $\mu, \nu$,
\[
H_{\mu, \nu} = \pi H[\nu] + \pi H[\mu] \frac{d\nu}{d \mu};
\]
note however that $H_{\mu, \nu}$ may be defined even if $H_\mu$, $H_\nu$ are not.
\end{Remark}

\begin{Lemma}
For any $\mu, \nu$ for which $H_{\mu, \nu}$ is well defined, the operator
\[
L_\nu L_\mu - H_{\mu, \nu} L_\mu
\]
is symmetric with respect to the inner product induced by $\mu$.
\end{Lemma}

\begin{proof}
\[
\begin{split}
\mu[(L_\nu[L_\mu[f]] - H_{\mu, \nu} L_\mu[f]) g]
& = (\mu \otimes \nu) [(\partial L_\mu[f]) (g \otimes I)] - \mu[H_{\mu, \nu} L_\mu[f] g] \\
& = - (\mu \otimes \nu) [(I \otimes L_\mu[f]) \partial g]
= - \nu[L_\mu[f] L_\mu[g]]. \qedhere
\end{split}
\]
\end{proof}

\begin{Remark}
\label{Remark:Phi}
The map $\mu = \Phi_{\beta, \gamma}[\nu]$ is the left shift of the Jacobi parameter sequences of the measure: it takes a measure with Jacobi parameters $\set{(\beta_0, \beta_1, \beta_2, \ldots,), (\gamma_1, \gamma_2, \gamma_3, \ldots)}$ to the measure with Jacobi parameters $\set{(\beta, \beta_0, \beta_1, \beta_2, \ldots,), (\gamma, \gamma_1, \gamma_2, \gamma_3, \ldots)}$, and so is an inverse of ``coefficient stripping'' \cite{Damanik-Simon-periodic}. Equivalently,
\[
M_\mu^{-1}(z) = 1 - \beta z - \gamma z^2 M_\nu(z)
\]
or
\[
G_\mu(z) = \frac{1}{z - \beta - \gamma G_\nu(z)}.
\]
The map $\Phi = \Phi_{0, 1}$  was defined in \cite{Belinschi-Nica-B_t}. In particular, the free Meixner distributions are exactly the image under $\Phi$ of general semicircular distributions. One can express $\Phi_{\beta, \gamma}[\nu]$ through $\Phi[\nu]$ using the operation of boolean convolution, see \cite{AnsEvolution}.
\end{Remark}

\begin{Lemma}
\label{Lemma:free-Appell}
A monic polynomial family $\set{A_n}$ such that
\[
L_\nu[A_n] = A_{n-1}
\]
is orthogonal with respect to some $\mu$ if and only if $\nu$ is a semicircular distribution and $\mu = \Phi_{\beta, \gamma}[\nu]$ is a free Meixner distribution.
\end{Lemma}

\begin{proof}
For any $\mu, \nu$, there is a unique monic polynomial family such that $L_\nu[A_n] = A_{n-1}$ and $\mu[A_n] = 0$ for $n \geq 1$. These are the c-free Appell polynomials $A^{\mu, \nu}_n$ investigated in \cite{AnsAppell3}, where the statement about their orthogonality was also proved.
\end{proof}

\begin{Lemma}
\label{Lemma:H=x}
$\mu = \Phi[\nu]$ if and only if $H_{\mu, \nu} = x$. More precisely, for compactly supported $\mu, \nu$, for any continuously differentiable $f$,
\[
(\nu \otimes \mu) \left[ \partial f \right] = \left[ x f \right]
\]
if and only if $\mu = \Phi[\nu]$. More generally, $H_{\mu, \nu} = \frac{1}{\gamma} (x - \beta)$ if and only if $\mu = \Phi_{\beta, \gamma}[\nu]$.
\end{Lemma}

\begin{proof}
This result was proved in \cite{AnsAppell3} for polynomial $f$. Here is a more direct analytic proof. If $\mu = \Phi[\nu]$, then
\[
G_\nu(z) = z - \frac{1}{G_\mu(z)}.
\]
Therefore
\begin{equation}
\label{Resolvent}
\begin{split}
(\nu \otimes \mu) \left[ \partial \frac{1}{z-x} \right]
& = \nu \left[ \frac{1}{z-x} \right] \mu \left[ \frac{1}{z-x} \right]
= \left( z - \frac{1}{G_\mu(z)} \right) G_\mu(z) \\
& = z G_\mu(z) - 1
= \mu \left[ \frac{z}{z-x} \right] - \mu[1]
= \mu \left[ x \frac{1}{z-x} \right].
\end{split}
\end{equation}
Since $\mu, \nu$ are compactly supported, this implies that for any $n$, $(\nu \otimes \mu) \left[ \partial x^n \right] = \mu \left[ x^{n+1} \right]$. By the Weierstrass theorem, $\mf{R}[x]$ is dense in $C^{(1)}[\supp{\mu}]$ with the uniform norm on $f'$. Note that if $P_n' \rightarrow f'$ uniformly on a compact set, by adjusting the constant term  we can arrange to have $P_n \rightarrow f$ uniformly. Now take $f$ continuously differentiable. The functional $f \mapsto (\nu \otimes \mu) \left[ \partial f \right]$ is continuous with respect to the uniform norm on $f'$, and the functional $f \mapsto \mu \left[ x f \right]$ is continuous with respect to the uniform norm on $f$. Therefore they coincide on any continuously differentiable function. The converse follows by running equation~\eqref{Resolvent} backwards, and the proof of the more general statement is similar.
\end{proof}

\begin{Prop}
\label{Prop:SL-two-measures}
For polynomial $p, q$, the operator
\[
Q = p(x) L_\nu L_\mu + q(x) L_\mu,
\]
has a polynomial system of eigenfunctions if and only if
\[
(a z^2 + b z + c) M_\mu(z) M_\nu(z) + (d z + e) M_\mu(z) = c + e,
\]
where $p(x) = a + b x + c x^2$, $q(x) = d + e x$, and $M_\mu(z)$, $M_\nu(z)$ are the moment generating functions of $\mu, \nu$, respectively.
\end{Prop}

\br
The proof is similar to that of Theorem~\ref{Thm:SL-one-measure}.

\begin{Remark}
The conclusion of Proposition~\ref{Prop:SL-two-measures} is equivalent to
\[
(c + e) \mc{R}^{\mu, \nu}(u) = c \mc{R}^\nu(u) - d u - a u^2 - b u \mc{R}^\nu(u) - c \mc{R}^\nu(u)^2,
\]
where $\mc{R}^{\nu}$ is the combinatorial $R$-transform of $\nu$, and $\mc{R}^{\mu, \nu}$ is the $R$-transform in the two-state free probability theory \cite{AnsAppell3}.
\end{Remark}

\begin{Example}
If $b = c = 0$, and taking without loss of generality $e = -1$, we get
\[
M_\mu(z)^{-1} = 1 - d z - a z^2 M_\nu(z),
\]
in other words $\mu = \Phi_{d, a}[\nu]$. Recall that by Lemma~\ref{Lemma:H=x} these are exactly the cases when $H_{\mu, \nu}$ is a polynomial of degree one. Thus for arbitrary $\nu$, one can choose a positive $\mu$ so that $Q$ has polynomial eigenfunctions. Moreover, in this case the eigenfunctions are orthogonal. In particular, this is the case for $\nu$ a shifted semicircular distribution and $\mu$ the corresponding free Meixner distribution.
\end{Example}

\begin{Prop}
\label{Prop:All-orthogonal}
For $\mu = \Phi_{\beta, \gamma} [\nu]$, the monic orthogonal polynomials $P^\mu_n$ with respect to $\mu$ are eigenfunctions for the operator
\[
Q = L_\nu L_\mu - \gamma^{-1} (x - \beta) L_\mu
\]
with eigenvalue $- \gamma^{-1}$ for $n \geq 0$. A posteriori, this operator is simply
\[
Q[f] = \gamma^{-1} (- f + \mu[f]).
\]
\end{Prop}

\begin{proof}
It is a classical result about monic orthogonal polynomials of the first and second kind that
\[
L_\mu[P^\mu_n] = (I \otimes \mu)[\partial P^\mu_n] = P^\nu_{n-1}.
\]
So using Lemma~\ref{Lemma:H=x},
\[
\begin{split}
\mu[(L_\nu[L_\mu[P^\mu_n]] - H_{\mu, \nu} L_\mu[P^\mu_n]) P^\mu_k]
& = - \nu[L_\mu[P^\mu_n] L_\mu[P^\mu_k]]
= - \nu[P^\nu_{n-1} P^\nu_{k-1}] \\
& = - \delta_{nk} \norm{P^\nu_{n-1}}_\nu^2.
\end{split}
\]
But using notation from Remark~\ref{Remark:Phi},
\[
\norm{P^\mu_{n}}_\mu^2 = \gamma \gamma_1 \ldots \gamma_{n-1} = \gamma \norm{P^\nu_{n-1}}_\nu^2,
\]
so we conclude that
\[
\mu[(L_\nu[L_\mu[P^\mu_n]] - H_{\mu, \nu} L_\mu[P^\mu_n]) P^\mu_k] = - \gamma^{-1} \delta_{nk} \norm{P^\mu_{n}}_\mu^2
\]
and
\[
L_\nu[L_\mu[P^\mu_n]] - H_{\mu, \nu} L_\mu[P^\mu_n] = - \gamma^{-1} P^\mu_n
\]
for $n \geq 1$.
\end{proof}

\br
The following result should be well-known, but we have not found it in the literature. The particular cases involving Chebyshev polynomials of the first kind, and the free Meixner polynomials \cite[Theorem~4]{AnsMeixner}, are indeed well known. It follows from the preceding proposition, but a direct proof is also straightforward.

\begin{Cor}
Let $\mu$ be a measure with Jacobi parameters $\set{(\beta, \beta_0, \beta_1, \beta_2, \ldots,), (\gamma, \gamma_1, \gamma_2, \gamma_3, \ldots)}$, $\nu$ the ``once-stripped'' measure with Jacobi parameters $\set{(\beta_0, \beta_1, \beta_2, \ldots,), (\gamma_1, \gamma_2, \gamma_3, \ldots)}$, and $\tau$ the twice-stripped measure with Jacobi parameters $\set{(\beta_1, \beta_2, \ldots,), (\gamma_2, \gamma_3, \ldots)}$. Then the corresponding monic orthogonal polynomials satisfy
\[
x P^\nu_n(x) = P^\mu_{n+1}(x) + \beta P^\nu_n(x) + \gamma P^\tau_{n-1}(x).
\]
\end{Cor}

\begin{Example}
If $\nu$ is a Bernoulli distribution,
\[
\nu = \alpha_1 \delta_{\beta_1} + \alpha_2 \delta_{\beta_2}
\]
with $\alpha_1 + \alpha_2 = 1$, then
\[
M_\nu(z) = \frac{\alpha_1}{1 - \beta_1 z} + \frac{\alpha_2}{1 - \beta_2 z} = \frac{1 - (\alpha_1 \beta_2 + \alpha_2 \beta_1) z}{1 - (\beta_1 + \beta_2) z + \beta_1 \beta_2 z^2}.
\]
Now take $e = -1$, $d = \alpha_1 \beta_2 + \alpha_2 \beta_1$, $c = 1$, $b = - (\beta_1 + \beta_2)$, $a = \beta_1 \beta_2$, so that
\[
(a z^2 + b z + c) M_\nu(z) + (d z + e) = 0 = c + e.
\]
It follows from Proposition~\ref{Prop:SL-two-measures} that for Bernoulli $\nu$, for \emph{any} $\mu$ there is an operator of the form
\[
Q = (a + b x + x^2) L_\nu L_\mu - (x - d) L_\mu
\]
with polynomial eigenfunctions.
\end{Example}

\begin{Prop}
Among all operators $p L_\nu L_\mu + q L_\mu$ with polynomial eigenfunctions, only the ones in Proposition~\ref{Prop:All-orthogonal} have orthogonal eigenfunctions.
\end{Prop}

\begin{proof}
By re-scaling, we may assume that $\mu$ has mean zero and variance one. Denoting by $\set{m^\nu_k}$ the moments of $\nu$, it follows from Proposition~\ref{Prop:SL-two-measures} that
\begin{equation}
\label{Mean-zero-2}
b + c m^\nu_1 + d = 0
\end{equation}
and
\begin{equation}
\label{Variance-one-2}
a + b m^\nu_1 + c + c m^\nu_2 + e = 0.
\end{equation}
Also, using the notation from the proof of Theorem~\ref{Thm:SL-one-measure},
\[
c B(z) = c z (1 - M_\mu(z) M_\nu(z))
\]
and
\[
(c + e) C(z) = z M_\mu(z) (a z M_\nu(z) + d(1 - M_\nu(z))).
\]
Thus
\[
\mu[x + \alpha] = m^\mu_1 + \alpha = m^\mu_1 + (d/e) = 0
\]
implies $\alpha = 0$, so $d = 0$ and $b = - c m^\nu_1$. Also
\[
\mu[x^2 + \beta_2 x + \gamma_2] = m^\mu_2 + m^\mu_1 \beta_2 + \gamma_2 = m^\mu_2 - (c/c) m^\mu_1 (m^\mu_1 + m^\nu_1) + \frac{a - d m^\nu_1}{c+e} = 0,
\]
so $a = d m^\nu_1 - (c+e) = - (c+e)$. It follows from equation~\eqref{Variance-one-2} that $b m^\nu_1 + c m^\nu_2 = 0 = c \left(m^\nu_2 - (m^\nu_1)^2 \right)$, so unless $\nu$ is a delta measure, $c = 0$ and $b = 0$.

\br
If $\nu = \delta_\xi$ with $M_\nu(z) = \frac{1}{1 - \xi z}$, then using Proposition~\ref{Prop:SL-two-measures} and the calculations above,
\[
\left( \frac{-(c+e) z^2 - c \xi z + c}{1 - \xi z} + e \right) M_\mu(z) = c + e
\]
and
\[
M_\mu(z) = \frac{(c + e) (1 - \xi z)}{- (c + e) z^2 - c \xi z + c + e - e \xi z} = \frac{1 - \xi z}{1 - \xi z - z^2}.
\]
It follows that $\mu$ is a (Bernoulli) free Meixner distribution, $\mu = \mu_{\xi, -1} = \State{\delta_\xi} = \State{\nu}$.
\end{proof}


\def\cprime{$'$}
\providecommand{\bysame}{\leavevmode\hbox to3em{\hrulefill}\thinspace}
\providecommand{\MR}{\relax\ifhmode\unskip\space\fi MR }
\providecommand{\MRhref}[2]{%
  \href{http://www.ams.org/mathscinet-getitem?mr=#1}{#2}
}
\providecommand{\href}[2]{#2}

\end{document}